\documentclass{amsart}
\usepackage{amsmath,amssymb}
\usepackage{graphicx,tikz}
\newtheorem{theorem}{Theorem}[section]

\newtheorem{corollary}[theorem]{Corollary}

\newtheorem{lemma}[theorem]{Lemma}

\theoremstyle{definition}
\newtheorem{definition}[theorem]{Definition}

\theoremstyle{remark}
\newtheorem{remark}[theorem]{Remark}

\newcommand{\al}{\alpha}
\newcommand{\be}{\beta}
\newcommand{\de}{\delta}
\newcommand{\ep}{\varepsilon}

\newcommand{\ga}{\gamma}

\newcommand{\la}{\lambda}

\newcommand{\si}{\sigma}

\newcommand\vka\varkappa
\newcommand{\De}{\Delta}
\newcommand{\Ga}{\Gamma}
\newcommand{\La}{\Lambda}
\newcommand{\Si}{\Sigma}
\newcommand{\Om}{\Omega}

\newcommand\cK{\mathcal K}
\newcommand\tcK{\widetilde{\mathcal K}}

\newcommand{\tq}{\tilde{q}}

\def\CC{\mathbb{C}}

\def\RR{\mathbb{R}}

\def\ZZ{\mathbb{Z}}

\def\TT{\mathbb{T}}

\renewcommand\SS{\mathbb{S}}
\newcommand{\cE}{{\mathcal E}}

\newcommand{\cI}{{\mathcal I}}
\newcommand{\cJ}{{\mathcal J}}
\newcommand{\cL}{{\mathcal L}}

\newcommand{\cP}{{\mathcal P}}

\newcommand{\cV}{{\mathcal V}}
\newcommand{\cU}{{\mathcal U}}

\newcommand\tQ{\widetilde{Q}}
\newcommand\tcE{\widetilde{\cE}}

\newcommand{\id}{{\rm id}}
\newcommand{\ms}{\mspace{1mu}}
\newcommand\lan\langle
 \DeclareMathOperator\GL{GL} \DeclareMathOperator\orto{O}

\DeclareMathOperator\dist{dist}

\renewcommand\leq\leqslant
\renewcommand\geq\geqslant
\newlength{\intwidth}

\numberwithin{equation}{section}

\newcommand{\triple}[1]{{\left\vert\kern-0.25ex\left\vert\kern-0.25ex\left\vert #1
        \right\vert\kern-0.25ex\right\vert\kern-0.25ex\right\vert}}

\begin{document}

\title{Localization
  properties of high energy eigenfunctions on flat tori}

 \author{Alberto Enciso}
 \address{Instituto de Ciencias Matem\'aticas, Consejo Superior de
   Investigaciones Cient\'\i ficas, 28049 Madrid, Spain}
 \email{aenciso@icmat.es}

 \author{Alba Garc\'\i a-Ruiz}
\address{Instituto de Ciencias Matem\'aticas, Consejo Superior de
   Investigaciones Cient\'\i ficas, 28049 Madrid, Spain}
 \email{alba.garcia@icmat.es}

 \author{Daniel Peralta-Salas}
 \address{Instituto de Ciencias Matem\'aticas, Consejo Superior de
   Investigaciones Cient\'\i ficas, 28049 Madrid, Spain}
 \email{dperalta@icmat.es}

%
%
\begin{abstract}
We consider the question of when the Laplace eigenfunctions on an arbitrary flat torus $\TT_\Ga:=\RR^d/\Ga$ are flexible enough to approximate, over the natural length scale of order $1/\sqrt\la$ where $\la\gg1$ is the eigenvalue, an arbitary solution of the Helmholtz equation $\De h + h=0$ on~$\RR^d$. This problem is motivated by the fact that, by the asymptotics for the local Weyl law, ``approximate Laplace eigenfunctions'' do have this approximation property on any compact Riemannian manifold. What we find is that the answer depends solely on the arithmetic properties of the spectrum. Specifically, recall that the eigenvalues of~$\TT_\Ga$ are of the form $\la_k=Q_\Ga(k)$, where $Q_\Ga$ is a quadratic form and $k\in\ZZ^d$. Our main result is that the eigenfunctions of~$\TT_\Ga$ have the desired approximation property if and only~$Q_\Ga$ is a multiple of a quadratic form with integer coefficients. In particular, the set of lattices~$\Ga$ for which this approximation property holds has measure zero but includes all rational lattices. A consequence of this fact is that when $Q_\Ga$ is a multiple of a quadratic form with integer coefficients, Laplace eigenfunctions exhibit an extremely flexible behavior over scales of order $1/\sqrt\la$. In particular, there are eigenfunctions of arbitrarily high energy that exhibit nodal components diffeomorphic to any compact hypersurface of diameter $O(1/\sqrt\la)$.
\end{abstract}

\maketitle

\section{Introduction}

The analysis of high energy eigenfunctions on compact Riemannian
manifolds in general, and on flat tori in
particular~\cite{Dima,BR,Riviere,Lagace} is a central
topic in geometric analysis where techniques from a number of areas of
mathematics come into play.

In this paper we are interested in high energy eigenfunctions on flat tori of dimension $d\geq2$. We shall use the notation $\TT_\Ga:=\RR^d/\Ga$ for the flat torus defined by a full-rank lattice $\Ga\subset\RR^d$. It is well known that the set of full-rank lattices on~$\RR^d$ is $\cL_d:= \GL_d(\RR)/\GL_d(\ZZ)$, and that the set of flat tori can be naturally identified with $\cL_d/\orto(d)$, that is, with the set of full-rank lattices up to an orthogonal transformation. Laplace eigenfunctions of~$\TT_\Ga$ satisfy the equation
\[
\De u + \la u=0
\]
on~$\RR^d$ for some nonnegative constant~$\la$ and the periodicity condition
$$
u(x +\ga)=u(x)
$$
for all $x\in\RR^d$ and all $\ga\in\Ga$.

The local behavior of a high energy eigenfunction on a flat torus (and, actually, on any compact Riemannian manifold) over length scales of order $1/\sqrt\la$ is described by a solution to the Helmholtz equation,
\begin{equation}\label{Helmholtz}
\De h + h =0	\,.
\end{equation}
In our setting, this is an obvious consequence of the fact that the rescaled eigenfunction $v(x):= u(x/\sqrt\la)$ does satisfy the Helmholtz equation.

A well known partial converse is that, given any ball $B\subset\RR^d$, one can pick a sequence of {\em approximate}\/ Laplace eigenfunctions on the flat torus $\TT_\Ga$ of approximate frequency~$\la$ whose behavior on a ball of radius $1/\sqrt\la $ reproduces that of any fixed solution~$h$ to the Helmholtz equation~\eqref{Helmholtz} on~$B$ modulo a small error.
More precisely, by the well known asymptotics for the local Weyl law~\cite{Hormander}, given any Helmholtz solution~$h$, any positive integer~$l$ and any~$\ep>0$, for any $\eta>0$ and all large enough $\la$ one can construct a linear combination~$U$ of eigenfunctions with eigenvalues in the interval $I_\la:=[(1-\eta)\la,(1+\eta)\la]$ such that
\begin{equation}\label{approxU}
\big\|U(\,\cdot\,/\sqrt\la)-h\big\|_{C^l(B)}<\ep\,.
\end{equation}
For the benefit of the reader, we sketch this construction (which works on any compact Riemannian manifold) in the Appendix. It is worth mentioning that the size of the interval $I_\la$ can be made substantially smaller, see e.g.~\cite{CanHen} and references therein.

A slightly imprecise but very intriguing question is when one can replace approximate eigenfunctions by bone fide eigenfunctions in the estimate~\eqref{approxU}. Our objective in this paper is to answer this question completely in the class of flat tori. To frame the problem, let us make the following definition:
	
\begin{definition}\label{D.IL}
	The flat torus $\TT_\Ga$ has the {\em inverse localization property} if for any solution of the Helmholtz equation~$h$ on~$\RR^d$, any ball $B\subset\RR^d$ and any positive integer~$l$, there exists a sequence of Laplace eigenfunctions~$u_n$ on~$\TT_\Ga$ with eigenvalues~$\la_n\to\infty$ whose rescalings approximate~$h$ as
	\begin{equation}\label{defapprox}
	\lim_{n\to\infty} \big\| u_n(\,\cdot\,/\sqrt{\la_n}) - h\big\|_{C^l(B)}=0\,.
	\end{equation}
\end{definition}

Several comments are in order. Firstly, note that Definition~\ref{D.IL} can be trivially extended to arbitrary compact Riemannian manifolds using the exponential map $\exp_p$ at a point $p$ in the manifold, and that in this more general setting it makes sense to use different base points~$p_n$ for the exponential. However, as the isometry group of~$\TT_\Ga$ is transitive and Laplace eigenspaces are invariant under isometries, in the case of flat tori that we study in this paper there is no loss of generality in only considering scalings centered at a fixed point, which we take to be the origin.

Also, let us point out that the technical details of Definition~\ref{D.IL} are unimportant in the sense that the specific norm and the specific sets used in the estimate~\eqref{approxU} (which we have chosen, for concreteness, as a certain~$C^l$ norm on balls~$B$) and whether the function~$h$ satisfies the Helmholtz equation on the whole~$\RR^d$ or only in a neighborhood of~$B$ are actually irrelevant. Roughly speaking, the choice of a norm is irrelevant because, on scales of order~$1/\sqrt\la$, the eigenfunction satisfies an elliptic equation with bounded coefficients, so standard elliptic estimates can be used to pass from bounds in a weak norm (say $L^2(B')$) to bounds in a stronger norm (say $C^l(B)$). The fact that the choice of the sets and the domain of definition of the function~$h$ are not important either comes from the existence of global approximation theorems for the Helmholtz equation, which permit to approximate solutions to the Helmholtz equation on sets satisfying certain mild technical hypotheses by solutions on the whole space. We will get back to this point in Section~\ref{S.GAT}.

The inverse localization property was established in~\cite{proceedings} for Laplace eigenfunctions on the round sphere $\mathbb S^d$ and the standard flat torus $\mathbb T^d$, and in~\cite{JEMS} for the quantum harmonic oscillator in Euclidean space. In the context of Beltrami fields (i.e., eigenfunctions of the curl operator) an analogous inverse localization property also holds~\cite{ENS} on $\mathbb S^3$ and $\mathbb T^3$. These results were used to construct high-energy eigenfunctions with nodal sets of complicated topologies. We want to emphasize that there are fundamental differences between the inverse localization property for eigenfunctions on $\mathbb S^d$ and $\mathbb T^d$. While the former holds for any sequence of eigenvalues $\lambda_n\to\infty$, the latter needs to restrict to special sequences enjoying certain number theoretic properties (related to the so called Linnik's problem). Moreover, the inverse localization on $\mathbb S^d$ exhibits a concentration property that allows us to extend it to quotients of $\mathbb S^d$ by finite isometry groups (lens spaces)~\cite{proceedings}; an analogous concentration result is not true on $\mathbb T^d$, which is consistent with the fact that the inverse localization property changes drastically if one considers flat tori $\TT_\Ga$ for different lattices.

Let us now state our results. For this, we record here that the Laplace eigenvalues of the flat torus~$\TT_\Ga$ are $\{4\pi^2|\xi|^2: \xi\in\Ga^*\}$, where $\Ga^*$ is the dual lattice. Therefore, the eigenvalues are labeled in terms of an integer vector~$k\in\ZZ^d$ as
\[
\la_k= Q_\Ga(k)\,,
\]
where $Q_\Ga$ is a positive-definite quadratic form defined by the lattice. The multiplicity of~$\la_k$ is given by the number of integral points on the ellipsoid $Q_\Ga^{-1}({\la_k})$. It can be proved that the multiplicity of the eigenvalues is unbounded if, but not only if, $Q_\Ga$ is a multiple of a quadratic form with integer coefficients. Recall that a quadratic form
\[
B(k)=\sum_{1\leq i\leq j\leq d} b_{ij}k_ik_j
\]
is said to have integer coefficients if $b_{ij}\in\ZZ$; the key point to note here is that the sum ranges over $i\leq j$ so the factor in front of, say, $k_1k_2$ is not necessarily even.

Our main result asserts that the inverse localization property holds if and only if~$Q_\Ga$ has integer coefficients up to a multiplicative constant:

\begin{theorem}\label{T.main}
	The $d$-dimensional flat torus $\TT_\Ga$ has the inverse localization property if and only if $Q_\Ga$ is a multiple of a quadratic form with integer coefficients.
\end{theorem}

In terms of the lattices~$\Ga$, one can analyze the validity of the inverse localization principle as follows:

\begin{theorem}\label{T.main2}
	For a full measure set of lattices, $\TT_\Ga$ does not have the inverse localization property. However, $\TT_\Ga$ does have the inverse localization property when a scaling of $\Ga$ is an integer lattice, that is, when $r\ms \Ga\subset\ZZ^d$ for some~$r>0$.
\end{theorem}

It should be noted that being able to approximate any solution of the Helmholtz equation implies that Laplace eigenfunctions exhibit an extremely flexible behavior over scales of order $1/\sqrt\la$. For example, in view of the results about solutions to the Helmholtz equation proved in~\cite{EPS13,CS19}, a manifold with the inverse localization property must have eigenfunctions of arbitrarily high energy that present a nodal component diffeomorphic to any given compact hypersurface of~$\RR^d$, and any fixed number of  nondegenerate critical points in a certain geodesic ball of radius $C/\sqrt\la$. We will make this assertion precise in Theorem~\ref{T.nodalsets}.

The paper is organized as follows. In the next two short sections, we shall recall some formulas for eigenfunctions on flat tori and a global approximation result for the Helmholtz equation that we will need later on. In Section~\ref{S.noIL} we show that if $Q_\Ga$ is not a multiple of a quadratic form with integer coefficients, then $\TT_\Ga$ cannot have the inverse localization property. The proofs of the converse implication and of Theorem~\ref{T.main2} are respectively presented in Sections~\ref{S.IL} and~\ref{S.C}. An application to the study of nodal sets and critical points of eigenfunctions on flat tori when $Q_\Ga$ has integer coefficients up to a constant is detailed in Section~\ref{S.applications}.

\section{Spectral theory of flat tori}
\label{S.flattori}

In this section we recall some well known facts about Laplace eigenfunctions on the $d$-dimensional flat torus $\TT_\Ga$.

Let us start by taking a basis $\{a_\Ga^j\}_{j=1}^d$ of the lattice~$\Ga$ and letting $A_\Ga\in \GL_d(\RR^d)$ be the matrix whose rows are the vectors $a_\Ga^j$. Equivalently, $A_\Ga$ is any matrix such that $\Ga= A_\Ga\ZZ^d$. The dual lattice,
\[
\Ga^*:=\{\xi\in\RR^d: \xi\cdot x\in\ZZ\; \text{for all } x\in\Ga\}\,,
\]
is then $\Ga^*=B_{\Ga}\ZZ^d$ with
\begin{equation}\label{AB}
	B_{\Ga}:=(A_\Ga^*)^{-1}\,.
\end{equation}
Here and in what follows, the star denotes the transposed matrix.

It is easy to see that the functions
\[
\{ e^{2\pi i x\cdot\xi} : \xi\in \Ga^*\} = \{ e^{2\pi i x\cdot B_\Ga k} :k\in\ZZ^d\}
\]
are an orthogonal basis of $L^2(\TT_\Ga)$. Since
\[
\De e^{2\pi i x\cdot\xi}= -4\pi^2 |\xi|^2 e^{2\pi i x\cdot\xi}\,,
\]
one infers that the spectrum of~$\TT_\Ga$, counting multiplicities, is $\{\la_k:k\in\ZZ^d\}$ with
\begin{equation}\label{Eqlak}
\la_k:=  Q_\Ga(k)\,.
\end{equation}
With some abuse of notation, let us denote by $Q_\Ga$ the symmetric matrix defined by the quadratic form via
\begin{equation}\label{matrixQ}
Q_\Ga(k):=k\cdot Q_\Ga k\,.
\end{equation}
We then have
\begin{equation}\label{QB}
Q_\Ga:= 4\pi^2B_\Ga^* B_\Ga\,,
\end{equation}
so that
\begin{equation}\label{QBGa}
Q_\Ga(k)=4\pi^2|B_\Ga k|^2\,.	
\end{equation}
Therefore, if one sets
\begin{equation}\label{mult}
\cK_k:=\{k'\in\ZZ^d: Q_\Ga(k')=Q_\Ga(k)\}\,,
\end{equation}
it follows that
\[
\{ e^{2\pi i x\cdot B_\Ga k'} :k'\in \cK_k\}
\]
is an orthogonal basis of the eigenspace~$\mathcal U_{\la_k}$ corresponding to the eigenvalue~$\la_k$.

\section{Helmholtz solutions and Fourier transforms of distributions supported on the unit sphere}
\label{S.GAT}

In some parts of the paper, it will be very useful to approximate solutions of the Helmholtz equation on a ball by global solutions of this equation that are given by the Fourier transform of an absolutely continuous measure on the unit sphere
\[
{\SS^{d-1}}:=\{\xi\in\RR^d:|\xi|=1\}
\]
with a nice density.

Let us recall that, by a classical result of Herglotz~\cite[Theorem 7.1.28]{Hor15}, the integrability properties of this density function are closely related to the decay at infinity of the solution. Specifically, any solution of the Helmholtz equation on~$\RR^d$ which satisfies the sharp $L^2$~decay condition
\begin{equation}\label{decayL2}
\triple{h}^2:=\limsup_{R\to\infty}\frac1R\int_{|x|<R}h(x)^2\, dx<\infty
\end{equation}
can be written as the Fourier transform of a measure of the form $H\, d\si$ with $H\in L^2({\SS^{d-1}})$, and that in fact the $L^2$~norm of~$H$ and the above seminorm are equivalent in the sense that
\[
\frac{\|H\|_{L^2({\SS^{d-1}})}}C\leq \triple{u} \leq C\|H\|_{L^2({\SS^{d-1}})}\,.
\]

The basic result about Helmholtz solutions of the form $h=\widehat{H\, d\si}$ that we will employ in this paper is the following, which was established in a slightly different form in~\cite{Acta,APDE}.
Throughout, $d\si$ will denote the spherical measure on~$\SS^{d-1}$.

\begin{theorem}\label{T.GAT}
Let $h_0$ satisfy the Helmholtz equation~\eqref{Helmholtz} in a neighborhood of the closure of a bounded open set~$\Om\subset\RR^d$ whose complement $\RR^d\backslash \Om$ is connected. Then for any positive integer~$l$ and any~$\ep>0$, there exists a complex-valued Hermitian polynomial~$p$ (i.e., $p(-\xi)=\overline{p(\xi)}$) such that
\[
h(x):=\int_{\SS^{d-1}} e^{i x\cdot\xi} p(\xi) \, d\si(\xi)
\]
approximates~$h_0$ as
\[
\|h-h_0\|_{C^l(\Om)}<\ep\,.
\]
Furthermore, $h$ satisfies the Helmholtz equation on~$\RR^d$ and the sharp pointwise decay condition
\begin{equation}\label{decay}
|h(x)|\lesssim
 (1+|x|)^{(1-d)/2}
\end{equation}
\end{theorem}

\begin{proof}
It was proved in~\cite{APDE} that this global approximation property holds for a function~$h$ given by a finite Fourier--Bessel series of the form
\[
h(x)= \sum_{l=0}^L \sum_{m=1}^{d_l} a_{lm}\,\frac{J_{l+\frac d2-1}(|x|)}{|x|^{\frac d2-1}}\, Y_{lm}\left(\frac x{|x|}\right)\,,
\]
for some real coefficients $a_{lm}$. Here $L$ is an integer, $J_\nu$ denotes the Bessel function of order~$\nu$ and
$\{Y_{lm}(\xi)\}$ is a real-valued orthonormal basis of $d$-dimensional spherical harmonics; the order~$l$ means that the spherical harmonic is the restriction to the sphere~$\SS^{d-1}$ of a homogeneous harmonic polynomial of degree~$l$, and $d_l$ is the multiplicity of this space.

Since
\[
b_l\,\widehat{Y_{lm}\, d\si}(x)=  \frac{J_{l+\frac d2-1}(|x|)}{|x|^{\frac d2-1}}\, Y_{lm}\left(\frac x{|x|}\right)
\]
for some explicit nonzero constant $b_l$ that does not depend on $m$ (see e.g.~\cite[Proposition 2.1]{Alvaro}), one finds that in fact one can write
\[
h(x)= \int_{\SS^{d-1}}\sum_{l=0}^L\sum_{m=1}^{d_l} a_{lm}b_l Y_{lm}(\xi) \, e^{ix\cdot\xi}\, d\si(\xi)=: \int_{\SS^{d-1}}p(\xi)\, e^{ix\cdot\xi}\, d\si(\xi)
\]
for some explicit constants $b_l$. Obviously, $p(\xi)$ is Hermitian because $h$ is real-valued.
\end{proof}

\begin{remark}\label{R.appendix}
The function~$h$ can be written as $h=P\cJ$, where $P:= p(-i\nabla)$ is a differential operator with constant coefficients and $\cJ$ is
the spherical Bessel function
\begin{equation}\label{defJ}
\cJ(x):= \int_{{\SS^{d-1}}} e^{i x\cdot \xi}\, d\si(\xi)= c_d |x|^{\frac12-d} J_{\frac d2-1}(|x|)\,,
\end{equation}
which is the only spherically symmetric solution to the Helmholtz equation on~$\RR^d$ up to a multiplicative factor. In this formula, $c_d$ is an explicit dimensional constant. Note that the pointwise decay condition~\eqref{decay} obviously implies the $L^2$-decay condition~\eqref{decayL2}.
\end{remark}

\section{Non-integer forms and failure of the inverse localization property}
\label{S.noIL}

Our objective in this section is to prove the following theorem, which asserts that flat tori for which $Q_\Ga$ is not a multiple of a quadratic form with integer coefficients do not possess the inverse localization property. This is a slightly more precise statement of the ``only if'' part of Theorem~\ref{T.main}.

\begin{theorem}\label{T.noIL}
	Fix a ball $B\subset\RR^d$. If $Q_\Ga$ is not a multiple of a quadratic form with integer coefficients, there is a solution of the Helmholtz equation~$h$ on~$\RR^d$ and some $\ep>0$ such that any Laplace eigenfunction~$u$ on~$\TT_\Ga$ satisfies
	\[
\int_B\big[h(x)-u(x/\sqrt\la)\big]^2\, dx >\ep\,.
	\]
	Here $\la$ denotes the eigenvalue of~$u$.
	\end{theorem}

\subsection{Proof of the theorem}

Next we present the proof of Theorem~\ref{T.noIL}. Consequently, in the rest of this section we shall assume that $Q_\Ga$ is not a multiple of a quadratic form with integer coefficients. For concreteness, we also assume without any loss of generality that the set~$B$ contains the origin of the coordinate system.
Throughout, we shall use the notation
\[
v(x):=u(x/\sqrt\la)
\]
for the rescaled eigenfunctions.

As we mentioned in the Introduction, one should observe that the choice of the norm in which one approximates the Helmholtz solution~$h$ is immaterial. This is because the difference $h-v$ obviously satisfies the Helmholtz equation on~$\RR^d$, so for any fixed~$l$ one has
\[
\|h-v\|_{C^l(B)}\leq C \|h-v\|_{L^2(B')}
\]
by standard elliptic estimates, provided that the closure of $B$ is contained in the interior of~$B'$.

By Equation~\eqref{Eqlak}, the Laplace eigenvalues of~$\TT_\Ga$ can be written as
\[
\la=  Q_\Ga(k)
\]
with $k\in\ZZ^d$, and the corresponding eigenspace consists of the real-valued linear combinations of the form
\[
u(x)=\sum_{k'\in \cK_k} \al_{k'} e^{2\pi i B_\Ga k'\cdot x}\,.
\]
with $\al_{k'}\in\CC$
and $\cK_k$ given by~\eqref{mult}.
By Equation~\eqref{QBGa}, the rescaled eigenfunctions are therefore the real-valued linear combinations
\begin{equation}\label{E.formv}
v(x)=\sum_{k'\in\cK_k} \al_{k'} e^{ix\cdot \xi_{k'}}
\end{equation}
with
\begin{equation}\label{eta}
\xi_{k'}:=\frac{B_\Ga k'}{|B_\Ga k'|}\in\SS^{d-1}\,.
\end{equation}

To further analyze the rescaled eigenfunction~$v$, we need to use the hypothesis that the form coefficients $q_{ij}$, defined as
\[
Q_\Ga(k)=:\sum_{1\leq i\leq j\leq d} q_{ij} k_i k_j\,,
\]
cannot be written as $q_{ij}= \be b_{ij}$ with $b_{ij}\in\ZZ$ and $\be\in \RR$. An efficient way of doing so is by noticing that this means that there is some integer $m\in[2,\frac12d(d+1)]$ for which one can write
\begin{equation}\label{formulaQ}
Q_\Ga(k)=\sum_{r=1}^m \be_r Q_r(k)\,,
\end{equation}
where the real numbers $\be_1,\dots,\be_m$ are independent over the integers (that is, the only solution to the equation
\[
\sum_{r=1}^m n_r\be_r=0
\]
with $(n_1,\cdots,n_m)=:n\in\ZZ^m$ is $n=0$), and the quadratic forms
\[
Q_r(k):=\sum_{(i,j)\in \cI_r} b_{ij}k_ik_j
\]
have integer coefficients (that is, $b_{ij}\in\ZZ$). Here $\cI_1,\dots, \cI_m$ is a partition of the set
\begin{equation}\label{cI}
	\cI:=\{(i,j): q_{ij}\neq0,\; 1\leq i\leq j\leq d\}
\end{equation}
(that is, the sets $\cI_r$ are disjoint and their union is~$\cI$). Since $Q_\Ga$ is positive definite, $d\leq \#\cI\leq \frac12d(d+1)$.

With these definitions in hand, the following result is obvious:

\begin{lemma}\label{L.Qr}
With $Q_\Ga$ as in~\eqref{formulaQ}, $k'\in \cK_k$ if and only if $Q_r(k')=Q_r(k)$ for all $1\leq r\leq m$.	
\end{lemma}

\begin{proof}
Since the equation $Q_\Ga(k')=Q_\Ga(k)$ can be written as
\[
\sum_{r=1}^m \be_r \big[Q_r(k')-Q_r(k)]=0
\]
and $Q_r(k')\in\ZZ$ for all $k'\in\ZZ^d$, the lemma follows from the independence of $\be_r$ over the integers.
\end{proof}

Theorem~\ref{T.noIL} then follows from the following lemma. In the statement we use the spherically symmetric solution of the Helmholtz equation $\cJ$ introduced in Remark~\ref{R.appendix}.

\begin{lemma}\label{L.J}
	With~$Q_\Ga$ as above, there is some~$\ep>0$ such that
	\[
	\|\cJ-u(\cdot /\sqrt\la)\|_{C^2(B)}>\ep
\]
for any eigenvalue~$\la$ and any eigenfunction $u\in\cU_\la$.
\end{lemma}

\begin{proof}
Suppose that there is a sequence of rescaled eigenfunctions $v_n$ as above, corresponding to eigenvalues $\la_n:= Q_\Ga(k_n)$ with $k_n\in\ZZ^d$, which approximate~$\cJ$ so that
\begin{equation}\label{error}
r_n(x):= \cJ(x)-v_n(x)\to 0
\end{equation}
in $C^2(B)$ as $n\to\infty$. Lemma~\ref{L.Qr} implies that $\cK_{k_n}=\cK^n$, with
\[
\cK^n:=\{ k\in\ZZ^d: Q_\Ga(k)=Q_\Ga(k_n),\; Q_1(k)=Q_1(k_n)\}\,.
\]
We have already seen that $v_n$ must be of the form
\[
v_n(x)= \sum_{k\in\cK^n} \al_k e^{i\xi_k\cdot x}
\]
for some complex constants~$\al_k$. As $\xi_k= 2\pi B_\Ga k/\sqrt{Q_\Ga(k)}$, it follows that
\begin{equation}\label{Peta}
P(\xi_k)= \frac{Q_1(k_n)}{Q_\Ga(k_n)}=:c_n
\end{equation}
for all $k\in\cK^n$, where $P$ is the quadratic form
\[
P(\xi):= \frac1{4\pi^2}\, Q_1(B_\Ga^{-1}\xi)\,.
\]
Also, note that the constants~$c_n$ are uniformly bounded as
\begin{equation}\label{boundcn}
|c_n|\leq \max_{\xi\in\SS^{d-1}}\frac{Q_1(\xi)}{Q_\Ga(\xi)}<\infty
\end{equation}
because the quadratic form $Q_\Ga$ is positive definite.

Let us denote by~$P$ and~$Q_1$, with some abuse of notation, the symmetric matrices defined by the corresponding quadratic forms via the identity~\eqref{matrixQ}. One then has
\[
P= \frac{1}{4\pi^2}(B_\Ga^{-1})^* Q_1B_\Ga^{-1}\,.
\]
This implies that $P$ cannot be a multiple of the identity, for if $P=\ga I_d$ with $\ga\in\RR$, then
\[
Q_1=4\pi^2\ga B_\Ga^* B_\Ga = \ga\ms Q_\Ga\,,
\]
which is absurd because the quadratic form $Q_\Ga$ is not a multiple of~$Q_1$ by~\eqref{formulaQ}.

Consider the second order differential operator with constant coefficients
\[
L:=P(-i\nabla )
\]
associated with this quadratic form. Equation~\eqref{Peta} implies that
\[
L v_n=\sum_{k\in\cK^n} \al_k P(\xi_k)\, e^{i\xi_k\cdot x}= c_n v_n\,.
\]
It is easy to check that the fact that the symmetric matrix~$P$ is not a multiple of the identity then implies that~$L\cJ$ (or the action of~$L$ on any nonconstant radial function) is not spherically symmetric.

Now, let us consider the action of the operator $L-c_n$ on the identity~\eqref{error}. One obtains
\[
(L-c_n)r_n= (L-c_n)\cJ\,.
\]
The left hand side is bounded on~$B$ as
\[
|(L-c_n)r_n|\leq C\|r_n\|_{C^2(B)}
\]
with $C$ independent of~$n$ by~\eqref{boundcn}, and thus tends uniformly to~0. However, $L\cJ$ is a fixed non-radial function that does not vanish on any open set (because it is analytic), while~$\cJ$ is a radial function. Therefore $L\cJ - c_n J$ cannot tend  to zero on~$B$, contradicting the hypothesis. The lemma is then proven.
\end{proof}

\begin{remark}
The choice of~$\cJ$ as the Helmholtz solution one tries to approximate is not incidental, since approximating~$\cJ$ is a sort of acid test for the inverse localization property. This is because, as we showed in Remark~\ref{R.appendix}, any Helmholtz solution can be approximated by a function of the form $P\cJ$, where $P$ is a differential operator with constant coefficients.
\end{remark}

\begin{remark}
A heuristic way (albeit non rigorous) of understanding the proof of Lemma~\ref{L.J} is the following. For all $k\in\cK^n$ the points $\xi_k$ lie on the intersection of $\mathbb S^{d-1}$ with the quadric $\{P(\xi)=c_n\}$; let us denote this set by $\cP_n$. Since the matrix associated to $P$ is not a multiple of the identity, $\cP_n$ is a set on $\mathbb S^{d-1}$ of codimension at least $1$. Taking into account that the function $\cJ$ is the Fourier transform of the Lebesgue measure on $\mathbb S^{d-1}$ it is reasonable to think that $\cJ$ cannot be approximated (on compact sets) by the Fourier transform of measures supported on $\cP_n$. The proof of Lemma~\ref{L.J} is a way to make this intuition rigorous.
\end{remark}

\subsection{The case of surfaces}

It is worth noticing that, when $d=2$, one can find a slightly easier proof of Theorem~\ref{T.noIL} that in fact provides some more refined information. In this subsection we shall explore this fact, always under the assumption that~$Q_\Ga$ is not a multiple of a form with integer coefficients.

The key observation, a variation of which will be key for the proof of Theorem~\ref{T.main2} as well, is the following:

\begin{lemma}\label{L.d2}
If $d=2$, all the eigenspaces of the Laplacian are of multiplicity at most~$4$.
\end{lemma}

\begin{proof}
We start by observing that, when $d=2$, the set $\cI$ defined in~\eqref{cI} consists of two or three elements and~$m$ is either~2 or~3. In either case, at least one of the sets~$\cI_r$, say~$\cI_1$, has exactly one element. Relabeling the coordinates and redefining the constant $\be_1$ if necessary, one can assume that
\[
Q_1(k)=k_1^2\qquad \text{or}\qquad Q_1(k)=k_1k_2\,.
\]

Let us first assume that $Q_1(k)=k_1^2$. By Lemma~\ref{L.Qr}, this ensures that any $k'\in\cK_k$ must satisfy $k_1'=\pm k_1$. If $m=3$, either $Q_2(k)$ or~$Q_3(k)$ is simply $k_2^2$, so one infers that $k_2'=\pm k_2$. Therefore, $\# \cK_k\leq 4$.

If $m=2$, the vectors $k'\in \cK_k$ must satisfy an equation of the form
\[
Q_2(k')= b_{22}k_2'^2 + b_{12}k_1'k_2'= Q_2(k)\,.
\]
For each value of~$k_1'=\pm k_1$, this is a quadratic equation for~$k_2'$, so it has at most two solutions. Therefore, $\#\cK_k\leq 4$ also in this case.

Let us now consider the remaining case, $Q_1(k)=k_1k_2$. If there are no forms~$Q_r$ of the form~$k_1^2$, this ensures that $m=2$ and
\[
Q_2(k)=b_{11}k_1^2+b_{22}k_2^2\,.
\]
Consider now a vector $k'\in\cK_k$. One can then multiply the equation $Q_2(k')=Q_2(k)$ by $k_1'^2$ and use the fact that $Q_1(k')=Q_1(k)$ to obtain a quartic equation for~$k_1'$:
\[
b_{11}k_1'^4+b_{22}Q_1(k)^2= Q_2(k) k_1'^2\,.
\]
This equation for $k_1'$ has at most $4$ solutions. As $k_1'$ determines $k_2'$, the lemma is proven.
\end{proof}

To show that the multiplicity bound $\#\cK_k\leq 4$ implies that the inverse localization property does not hold on~$\TT_\Ga$ we  use the following fact:

\begin{lemma}\label{L.mult}
If $\# \cK_k\leq 4$ for all $k$, then for any open set $B\subset\RR^2$ there exists a solution~$h$ of the Helmholtz equation on the plane and some~$\ep>0$ such that
\[
\|h-u(\cdot/\sqrt\la )\|_{L^2(B)}>\ep
\]
for any eigenvalue~$\la$ of the $2$-dimensional torus $\TT_\Ga$ and any eigenfunction $u\in\cU_\la$.
\end{lemma}

\begin{proof}
	Let us define the set
	\[
\cV:= \left\{  \sum_{j=1}^4 \al_j\, e^{i z_j\cdot x}: \al_j\in\CC\,,\; z_j\in\RR^2\,,\; |z_j|=1\right\}\,,
	\]
	which can be understood as a (closed) submanifold of~$L^2(B)$ of real dimension $12$. It is clear that one can take linear combinations of plane waves on~$\RR^2$ with frequencies of unit norm whose distance to~$\cV$ (measured in $L^2(B)$) is positive. An example is
	\[
	h(x):= \sum_{j=1}^{J} \cos\left( x_1\cos\frac{2\pi j}J +x_2\sin\frac{2\pi j}J \right)
	\]
	with $J\geq5$.
		
	Note that any rescaled eigenfunction~$v$ corresponding to an eigenvalue~$\la_k$ belongs to~$\cV$ because~$v$ is of the form
	\[
	v(x)=\sum_{k'\in\cK_k} \al_{k'} e^{i\xi_{k'}\cdot x}
	\]
	and $\#\cK_k\leq 4$. Therefore, 	
	\[
	\inf_{k\in\ZZ^2} \inf_{u\in\mathcal U_{\la_k}}\left\|h-u(\cdot/\sqrt{\la_k})\right\|_{L^2(B)} \geq \dist_{L^2(B)}(h,\cV)>0\,.
	\]
	The lemma then follows.
\end{proof}

\begin{remark}\label{R.mult}
The proof remains valid in higher dimensions with trivial modifications. Therefore, if $\Ga$ is a $d$-dimensional lattice and $\#\cK_k\leq c$ for all~$k\in\ZZ^d$, then $\TT_\Ga$ does not have the inverse localization property.
\end{remark}

Note that the strategy of proof we have presented in the case of surfaces cannot be easily extended to higher dimensions. Indeed, for $d\geq3$ there are lattices, such as the one spanned by the vectors $\{2^{-1/4}e_1, e_2,\dots, e_d\}$ and associated with the quadratic form $Q_\Ga(k)= \sqrt2\,k_1^2+k_2^2+\cdots+ k_d^2$, for which the multiplicities of the eigenvalues are unbounded yet $\TT_\Ga$ does not have the inverse localization property.

\section{Integer forms and Linnik's problem}
\label{S.IL}

In this section we consider the case
\begin{equation}\label{QB2}
Q_\Ga= 4\pi^2 \be \,Q\,,	
\end{equation}
where $Q$ is a quadratic form with integer coefficients and $\be>0$. Our goal is to prove the following result, which corresponds to the ``if'' part of Theorem~\ref{T.main}.

\begin{theorem}\label{T.IL}
With $Q_\Ga$ as in~\eqref{QB2}, there are sequences of integers $n_j\to\infty$ such that $\la_j:=4\pi^2\be n_j$ is an eigenvalue and
\[
\lim_{j\to0} \inf_{u\in \cU_{\la_j}}\left\| h - u(\cdot/\sqrt{\la_j})\right\|_{C^l(B)}=0
\]
for any solution~$h$ of the Helmholtz equation on~$\RR^d$.
\end{theorem}

To prove this result, we start by noting that, as a consequence of Theorem~\ref{T.GAT}, for the purposes of Theorem~\ref{T.IL} there is no loss of generality in assuming that
\[
h(x)=\int_{\SS^{d-1}}e^{ix\cdot\xi} p(\xi)\, d\si(\xi)\,,
\]
where $p$ is a polynomial.

To prove Theorem~\ref{T.IL} we will use equidistribution properties of integer points on ellipsoids. We shall say that
\[
\cK^n:=\{k\in\ZZ^d: Q(k)=n\}
\]
is {\em asymptotically equidistributed}\/ over the ellipsoid
\begin{equation}\label{defcE}
\cE:=\{x\in\RR^d: Q(x)=1\}	
\end{equation}
along a certain sequence of integers $n_j$ if, for every function $f\in C^\infty(\RR^d)$,
\[
\lim_{j\to\infty}\frac1{\#\cK^{n_j}}\sum_{k\in \cK^{n_j}} f(k/\sqrt {n_j}) = \frac1{|\cE|} \int_{\cE} f(\eta)\, d\si_\cE(\eta)\,,
\]
where $d\si_\cE$ is the hypersurface measure.

The way we will use equidistribution is through the following lemma. In the statement, the Jacobian function $\cJ_\cE$ is defined through the transformation formula for measures
\[
d\si(\xi) = J_\cE\, d\si_\cE(\eta)
\]
when $\xi=:B\eta$, $B:= \be^{-1/2} B_\Ga$ and~$\xi\in\SS^{d-1}$. Note that, in this case,
\[
1=|\xi|^2= Q(\eta)\,.
\]
Since the transformation between $\xi$ and $\eta$ is linear, it is clear that the Jacobian function is simply the constant $\cJ_\cE=|\det B|\neq 0$.
Also, in the statement we make use of the unit vectors defined for each $k\in\cK^n$ as in~\eqref{eta}, which one can write as
\[
\xi_k:= \frac{B k}{\sqrt n}\,.
\]

\begin{lemma}\label{L.approx}
	If $\cK^n$ is asymptotically equidistributed over the ellipsoid~\eqref{defcE} along a certain sequence of integers, then
	\begin{equation}\label{conv}
	\frac{|\det B|\cdot|\cE|}{\# \cK^n} \sum_{k\in\cK^n} p(\xi_k)\, e^{ix\cdot\xi_k} \to \int_{\SS^{d-1}}e^{ix\cdot\xi} p(\xi)\, d\si(\xi)
	\end{equation}
	in $C^l(B)$ as $n\to\infty$ along that sequence. Here $p$ is any polynomial.
\end{lemma}

\begin{proof}
For any multiindex $\alpha=(\alpha_1,\cdots,\alpha_d)$ of nonnegative integers, we denote by $(B\eta)^\alpha$ the polynomial given by $(B\eta)_1^{\alpha_1}\cdots (B\eta)_d^{\alpha_d}$, where $(B\eta)_j$ is the $j$-th component of the vector $B\eta$. Since
\[
\big\{(B\eta)^\al \, e^{ix\cdot B\eta} p(B\eta)\, J_\cE(\eta)|_{\cE}: x\in B,\; |\al|\leq l\big\}
\]
is a uniformly bounded subset of~$C^m(\cE)$ for any~$m$, the lemma follows directly from the identity
\[
\int_{\SS^{d-1}} e^{ix\cdot\xi} p(\xi)\, d\si(\xi)= \int_{\cE} e^{ix\cdot B\eta} p(B\eta)\, J_\cE(\eta)\, d\si_\cE(\eta)=|\det B|\int_{\cE} e^{ix\cdot B\eta} p(B\eta)\, d\si_\cE(\eta)
\]
and the equidistribution hypothesis.
\end{proof}

Since the left hand side of~\eqref{conv} is a rescaled eigenfunction with eigenvalue $4\pi^2\be n$ by~\eqref{E.formv}, to prove Theorem~\ref{T.IL} it suffices to show that the integer points on~$Q^{-1}(n)$ are asymptotically equidistributed as $n\to\infty$ along a certain sequence of integers. This is a problem of Linnik type in analytic number theory. For $d\geq3$ it is standard that there are sequences of integers $n$ with the desired equidistribution property. Specifically, for $d\geq4$ one can use the Hardy--Littlewood circle method to show~\cite[Theorem 11.5]{Iwaniec} that integer points on $Q^{-1}(n)$ are asymptotically equidistributed as $n\to\infty$ over integers satisfying the congruence
\begin{equation*}
n=Q(k_0)\mod 2^7(\det Q')^3	
\end{equation*}
for some $k_0\in\ZZ^d$, while for $d=3$ one must additionally impose~\cite[Theorem 11.6]{Iwaniec} that $n$ be square-free and that the greatest common divisor of~$n$ and $2|\det Q'|$ be~1. Here $Q'$ is the $d\times d$ symmetric matrix with integer components defined by the quadratic form via $Q(k)=:\frac12 k\cdot Q'k$. This proves Theorem~\ref{T.IL} when $d\geq3$.

The case $d=2$ is more involved. To tackle this case, we will employ a result of Cilleruelo and C\'ordoba~\cite{CC}  asserting that, if $q$ is a square-free positive integer, there exists a sequence of integers $N_j\to\infty$ along which the integer points on the two-dimensional ellipsoid
\[
K_1^2+ qK_2^2 = N_j
\]
are asymptotically equidistributed. More precisely (see~\cite[Section II.D]{CC}), the integers are of the form
\[
N_j:= \prod_{m=1}^{\lfloor j e^{4 \sqrt q}\rfloor} (qm^2+1)\,.
\]

To reduce our problem to this case, let us start by writing out the components of the quadratic form~$Q$,
\[
Q(k)= q_{11}k_1^2+ q_{12}k_1k_2+q_{22}k_2^2\,,
\]
with $q_{ij}\in\ZZ$. Let us consider the determinant
\[
\tq:= 4 q_{11}q_{22}-q_{12}^2\,.
\]
Since the quadratic form is positive, $\tq$ is a positive integer. Let us write $\tq= p^2 q$, where $q$ is a square-free integer and $p$ is a positive integer. Note that $p=1$ when $\tq$ is square-free. Let us now introduce the linear change of variables $k=:TK$ defined by
\begin{align*}
k_1=:	p K_1-q_{12}K_2\,,\qquad k_2=: 2 q_{11} K_2
\end{align*}
and note that $T$ is a matrix with integer components. A short computation reveals that
\begin{equation}\label{QtQ}
Q(TK)=p^2q_{11}\, \tQ(K)
\end{equation}
with
\[
\tQ(K):=K_1^2+q K_2^2\,.
\]

Consequently, if we write
\begin{equation}\label{defn}
n:= p^2q_{11} N
\end{equation}
and restrict our attention to those $k\in\ZZ^2$ that can be written as $k=TK$ with $K\in\ZZ^2$, the equation $Q(k)=n$ reduces to
\begin{equation}\label{CC}
\tQ(K)=N\,.
\end{equation}
Furthermore, one has the relation
\[
|BTK|^2=Q(TK)=p^2q_{11}\,\tQ(K)
\]
for all $K\in\RR^d$. We will also be interested in the matrix
\[
M:=\frac1{p\sqrt{q_{11}}}BT\,,
\]
which is related to the quadratic form~$\tQ$ through the identity
\begin{equation}\label{defM}
\tQ(K)=|MK|^2\,.
\end{equation}

As $q$ is square-free, the aforementioned result of Cilleruelo and C\'ordoba ensures that there is a sequence of integers~$N_j$ for which the integer solutions
\[
\tcK^{N_j}:=\big\{K\in \ZZ^2: \tQ(K)=N_j\big\}
\]
to Equation~\eqref{CC} are asymptotically equidistributed on the ellipse
\[
\tcE:=\big\{\eta\in\RR^2: \tQ(\eta)=1\big\}\,.
\]
We shall denote by $d\si_{\tcE}$ the corresponding measure and introduce the Jacobian $J_{\tcE}$ as above, i.e., through the formula
\[
d\si(\xi)=J_{\tcE}\, d\si_{\tcE}(\eta)
\]
when $\xi= M\eta\in\SS^1$. Note $\xi\in\SS^1$ if and only if~$\eta\in\tcE$ by~\eqref{defM}. As before, the Jacobian  is simply the constant
\[
J_{\tcE}=|\det M|=\frac{|\det B||\det T|}{p^2q_{11}}\,.
\]

As $T$ maps $\ZZ^2$ into~$\ZZ^2$, it is clear in view of~\eqref{QtQ} that
\begin{equation}\label{TcK}
	T\tcK^N\subset \cK^n
\end{equation}
with $n$ given by~\eqref{defn}. Since
\[
\frac{MK}{\sqrt N}=\frac{BTK}{|BTK|}=\xi_{TK}
\]
with $\xi_k$ defined by~\eqref{eta},
one can apply Lemma~\ref{L.approx} to the asymptotically equidistributed sets~$\tcK^N$ on the ellipsoid~$\tcE$ to conclude that
\begin{align*}
\frac{|\det M|\cdot|\tcE|}{\# \tcK^N} \sum_{K\in\tcK^N} p(\xi_{TK})\, e^{ix\cdot\xi_{TK}} &= \frac{|\det M|\cdot|\tcE|}{\# \tcK^N} \sum_{k\in T\tcK^N} p(\xi_k)\, e^{ix\cdot\xi_k} \\
&\to \int_{\SS^{1}}e^{ix\cdot\xi} p(\xi)\, d\si(\xi)
\end{align*}
in $C^l(B)$ as $N\to\infty$ along the sequence~$N_j$. By~\eqref{TcK}, the right hand side of the first line is a rescaled eigenfunction in~$\cU_{\be n}$, thereby completing the proof of Theorem~\ref{T.IL} when $d=2$.

\section{Proof of Theorem~\ref{T.main2}}
\label{S.C}

The main idea of the proof of Theorem~\ref{T.main2} is to study the case where the coefficients of the quadratic form~$Q_\Ga$ satisfy a Diophantine condition. Recall that a vector $z\in\RR^N$ is~{Diophantine} if there exist constants $c>0$ and $\tau>N-1$ such that
\[
|K\cdot z|>c|K|^{-\tau}
\]
for all~$K\in\ZZ^N\backslash\{0\}$. It is well known that the set of Diophantine vectors has full measure on~$\RR^N$.

With some abuse of notation, let us denote by $Q_\Ga=(q_{ij})$ the $d\times d$ symmetric matrix associated with the quadratic form $Q_\Ga$ via~\eqref{matrixQ}. Understanding~$Q_\Ga$ as a $\frac12d(d+1)$-component vector, we will say that the quadratic form~$Q_\Ga$
is {\em Diophantine}\/ if there are constants $c>0$ and~$\tau>\frac12d(d+1)-1$ such that
\[
\left|\sum_{1\leq i\leq j\leq d} q_{ij} K_{ij}\right| > c\left(\sum_{1\leq i\leq j\leq d} |K_{ij}|\right)^{-\tau}
\]
for all $(K_{ij})_{1\leq i\leq j\leq d}\in\ZZ^{d(d+1)/2}\backslash\{0\}$.

An important first observation is the following:

\begin{lemma}\label{L.Dioph}
If $Q_\Ga$ is Diophantine, the multiplicity of each eigenvalue is at most~$2^d$.	
\end{lemma}

\begin{proof}
In the notation of Equation~\eqref{formulaQ}, the assumption that $Q_\Ga$ is Diophantine ensures that $m=\frac12d(d+1)$. Therefore, Lemma~\ref{L.Qr} ensures that
\[
k_ik_j= k'_ik'_j
\]
for all $k'$ belonging to the set~$\cK_k$ defined in~\eqref{mult}. Thus $k'\in\cK_k$ if and only if $k'_j=\pm k_j$, so necessarily $\#\cK_k\leq 2^d$ for all $k\in\ZZ^d$.	
\end{proof}

\begin{corollary}\label{C.Dioph}
	If $Q_\Ga$ is Diophantine, $\TT_\Ga$ does not have the inverse localization property.
\end{corollary}

\begin{proof}
An immediate consequence of Lemma~\ref{L.Dioph} and Remark~\ref{R.mult}.
\end{proof}

In view of these results, to prove Theorem~\ref{T.main2} it suffices to show that a Diophantine condition is satisfied for almost all lattices:

\begin{lemma}\label{L.measure}
For a full measure set of lattices, $Q_\Ga$ is Diophantine.	
\end{lemma}

\begin{proof}
By Equation~\eqref{QB},
\[
q_{ij}= 4\pi^2\sum_{k=1}^d b_{ki}b_{kj}\,,
\]
where $B_\Ga=(b_{ij})$ is the matrix associated with the dual lattice~$\Ga^*$. Let us henceforth identify the space of $d\times d$ matrices (to which $B_\Ga$ belongs) with $\RR^{d^2}$ and denote a typical element of this space by $B=(b_I)_{I=1}^{d^2}$. Note that the map
\[
F:\RR^{d^2}\to \RR^{d^2(d^2+1)/2}
\]
defined by
\[
F(B):=( b_I b_J)_{1\leq I\leq J\leq d^2}
\]
is analytic and that its image $F(\RR^{d^2})$ cannot be contained in any proper linear subspace of $\RR^{d^2(d^2+1)/2}$ because~$F$ is a quadratic function whose components are linearly independent monomials. This is well known to imply~\cite{Sevryuk} that the vector~$F(B)$ is Diophantine for almost all~$B\in \RR^{d^2}$. Therefore, for a full measure set of dual lattices~$\Ga^*$, $Q_\Ga$ is Diophantine.

	To conclude, we notice that if a set of dual lattices has full measure, the corresponding set of lattices has full measure as well. This is because, by Equation~\eqref{AB}, the matrix~$A_\Ga$ defining the lattice can be written as
	\[
	A_\Ga:= G(B_\Ga)\,,
	\]
	with $G(B):=(B^*)^{-1}$  being a diffeomorphism $\RR^d\backslash Z\to \RR^d\backslash Z$. Here
	\[
	Z:=\{B\in\RR^{d^2}: \det B=0\}
	\]
	is the set of non-invertible matrices, which has measure zero.
	
	The lemma, and thus Theorem~\ref{T.main2}, are then proven.
\end{proof}

\section{Applications to nodal sets and critical points of eigenfunctions}
\label{S.applications}

In this section we shall present a simple application of the inverse localization property to the study of nodal sets and critical points of eigenfunctions on flat tori that we alluded to in the Introduction.

Let us begin by introducing some notation. If $X$ is a subset of~$\RR^d$, we denote its dilation by factor $c>0$ and its translation by $p\in\RR^d$ by
\[
cX:=\{cx: x\in X\}\,,\qquad X+ p:=\{x+p: x\in X\}\,.
\]
Also, here $B_R$ denotes the ball centered at the origin of radius~$R$, which we can understand as a subset of~$\RR^d$ or of~$\TT_\Ga$ when $R$ is small enough.

In the statement of the following theorem we say that a collection of pairwise disjoint compact hypersurfaces~$\{\Si_j\}_{j=1}^N\subset\mathbb R^d$ is not linked if there are pairwise disjoint contractible domains $V_j$ such that $\Sigma_j\subset V_j$.

\begin{theorem}\label{T.nodalsets}
If a multiple of~$Q_\Ga$ is a quadratic form with integer coefficients, there exists a sequence of eigenvalues $\la_n\to\infty$ such that, given any~$N$ and a collection of compact embedded hypersurfaces~$\Si_j$ ($1\leq j\leq N$) of~$\RR^d$ that are not linked, any positive integer~$l$ and any $\ep$, there exists some~$R>0$ such that for all large enough~$n$ there is an eigenfunction $u_n$ with eigenvalue~$\la_n$ having in the ball $B_{R\la_n^{-1/2}}$ at least~$N$ nodal components of the form
\[
\widetilde\Si_j^n:= \la_n^{-1/2} \,  \Phi_{n}(c_j \Si_j + p_j)
\]
and at least $N$ nondegenerate local extrema. Here $c_j>0$, $p_j\in\RR^d$, and $\Phi_n$ is a diffeomorphism of~$\RR^d$ which is close to the identity:
$  \|\Phi_n-\id\|_{C^l(\RR^d)}<\ep$.
\end{theorem}

\begin{proof}
Let us choose the constants $c_j$ so that the first Dirichlet eigenvalue of the domain~$\Om_j$ bounded by the rescaled hypersurface $c_j\Si_j$ is~1. Let us now pick any vectors $p_j\in\RR^d$ so that the translated domains $\Om_j':=\Om_j+p_j$ have disjoint closures (which is possible because $\Sigma_j$ are not linked), and take some large ball~$B_R$ so that $\Om_j'\subset B_R$ for all $j$. It was proved in~\cite{EPS13,APDE} that then there exists a solution~$h$ to the Helmholtz equation which has structurally stable nodal sets diffeomorphic to
\[
\Si_j':=c_j \Si_j+p_j
\]
and at least $N$ nondegenerate local extrema in~$B_R$. The reason why we can assume that the local extrema are nondegenerate is that for generic domains, the first Dirichlet eigenfunctions are Morse~\cite{Uhlen}. More precisely, there exists some~$\de>0$ such that any function~$v$ on~$\RR^d$ with
\[
\|v-h\|_{C^{l}(B_R)}<\de
\]
has at least $N$ nondegenerate local extrema in~$B_R$ and the zero level set $v^{-1}(0)$ has at least $N$ components of the form $\Phi(\Si_j')$, where $\Phi$ is a diffeomorphism of~$\RR^d$ with $\|\Phi-\id\|_{C^l(\RR^d)}<\ep$.

By Theorem~\ref{T.IL}, there exists a sequence of rescaled eigenfunctions $v_n:=u_n({\cdot}/\sqrt{\la_n})$ with eigenvalues $\la_n\to\infty$ such that
\[
\lim_{n\to\infty}\|v_n-h\|_{C^{l+1}(B_R)}=0\,,
\]
so the result then follows from the above stability property of the function~$h$.
\end{proof}

\section*{Acknowledgements}

This work has received funding from the European Research Council (ERC) under the European Union's Horizon 2020 research and innovation programme through the grant agreement~862342 (A.E.\ and A.G.-R.). D.P.-S.
is supported by the grants MTM PID2019-106715GB-C21 and Europa Excelencia EUR2019-103821 from Agencia Estatal de Investigaci\'on (AEI). This work is supported in part by the ICMAT--Severo Ochoa grant CEX2019-000904-S and the grant RED2018-102650-T funded by MCIN/AEI/10.13039/501100011033.

\appendix

\section{Approximation of Helmholtz solutions by approximate eigenfunctions}
\label{S.appendix}

For the benefit of the reader, we shall next recall the standard construction of approximate eigenfunctions with the property~\eqref{approxU}; for the probabilistic counterpart of this result see e.g.~\cite{SW}. By Theorem~\ref{T.GAT}, we can assume that the Helmholtz solution is of the form
\[
h(x)=\int_{\SS^{d-1}}p(\xi)\, e^{ix\cdot \xi}\,d\si(\xi)\,,
\]
with $p$ a polynomial.

In this Appendix, $(M,g)$ will be a compact Riemannian manifold of dimension~$d$. If $u_j$ is an orthonormal basis of Laplace eigenfunctions, satisfying the equation
\[
\De_M u_j+ \la_j u_j=0
\]
on~$M$, the {\em spectral function}\/ of~$M$ corresponding to the interval~$I\subset \RR$ is
\begin{equation}\label{defEI}
E_I(x,y):=\sum_{{\la_j}\in I}u_j(x) \, u_j(y)\,.
\end{equation}
Note this is simply the kernel of the spectral projector onto the direct sum of the eigenspaces corresponding to the eigenvalues $\la_j\in I$.

Consider now the interval $I_\La:=[(1-\eta)^2\La^2,(1+\eta)^2\La^2]$, with $\La\gg1$, and denote by $d_M(x,y)$ the distance between two points $x,y\in M$, which we assume that are close enough. Consider points of the form
\begin{equation}\label{xy}
x=\exp_{x_0}(X/\La)\,,\qquad y= \exp_{x_0}(Y/\La)\,,
\end{equation}
where $x_0$ is a fixed point in~$M$ and $X,Y$ lie on the ball~$B_R\subset\RR^d$ of radius~$R\gg1$. H\"ormander's local Weyl law~\cite{Hormander} ensures that, in this case,
\begin{align*}
E_{I_\La}(x,y)&=\left(\frac\La{2\pi}\right)^d \int_{1-\eta<|\xi|_{g_y}<1+\eta} e^{i\La g_y(\xi,\exp_y^{-1}x)} \frac{d\xi}{\sqrt{\det g_y}} + O_R(\La^{d-1})\\
&=\left(\frac\La{2\pi}\right)^d \int_{1-\eta<|\xi|<1+\eta} e^{i(X-Y)\cdot \xi} \,{d\xi} + O_R(\La^{d-1})\,.
\end{align*}
Notice that $|\exp_y^{-1}x|=d_M(x,y)$. The subscript emphasizes that the error is not uniform in~$R$. The reader can consult e.g.~\cite{CanHen} for a lucid  summary of more recent refinements, which we will ignore in our discussion.

Now let~$F(\xi)$ be any smooth function on~$\RR^d$ such that $F(\xi):= p(\xi/|\xi|)$ if $|\xi|\in (\frac12,2)$ and $F(\xi):=0$ if $|\xi|<\frac14$ or $|\xi|>4$. As its Fourier transform is Schwartz, it is then clear that
\[
\int_{|Y|>R} |\widehat{\xi^\al\, F}(Y)|\, dY\leq C_N R^{-N}
\]
for any $N>0$ and any multiindex~$\al$. With $K$ a constant to be specified later, let us now define the function
\[
f_\La(y):= K\,\widehat F(\La \exp_{x_0}^{-1}y)\, \chi_1\left(\frac \La R \exp_{x_0}^{-1}y\right)\,,
\]
where $\chi_1$ is the indicator function of the ball of radius~$1$. Note that the function $f_\La$ is supported on a geodesic ball centered at~$x_0\in M$ of radius~$R/\La$ and is possibly discontinuous across the boundary of this ball.

Denoting by $dV$ the Riemannian measure on~$M$, consider the function
\[
U(x):= \int_M E_{I_\La}(x,y)\, f_\La(y)\, dV(y)\,.
\]
This is a linear combination of eigenfunctions with energies in the interval~$I_\La$ as a consequence of the identity~\eqref{defEI}. Furthermore, it satisfies the asymptotic estimate
\begin{align*}
U(x)&=K\left(\frac\La{2\pi}\right)^d \int_{|Y|<R}\int_{1-\eta<|\xi|<1+\eta} \widehat F(Y)\,e^{i(X-Y)\cdot \xi} \,{d\xi}\, dY + K\eta\,O_R(\La^{d-1})\\
&=K\left(\frac\La{2\pi}\right)^d \left[\int_{\RR^d}\int_{1-\eta<|\xi|<1+\eta} \widehat F(Y)\,e^{i(X-Y)\cdot \xi} \,{d\xi}\, dY +O(R^{-N})\right]+ K\eta\,O_R(\La^{d-1})\,.
\end{align*}
Here we are still assuming that $x,y$ are of the form~\eqref{xy}.
With~$\eta\ll1$, one can use the Fourier inversion formula and the fact that~$F(\xi)=p(\xi/|\xi|)$ for $|\xi|$ near~1 to write
\begin{align*}
U(x)&=K\La^d \left[\int_{1-\eta<|\xi|<1+\eta}  p\left(\frac\xi{|\xi|}\right)\,e^{iX\cdot \xi} \,{d\xi} +O(R^{-N})\right]+ K\eta\,O_R(\La^{d-1})\\
&=2\eta K\La^d \left[h(X)+ O(\eta) +O(R^{-N})\right]+ K\eta\,O_R(\La^{d-1})\,.
\end{align*}
Therefore, setting
\[
K:= \frac1{2\eta \La^d}
\]
one concludes that one can take a small enough~$\eta$ and large enough~$R$ and~$\La$ so that
\[
\left\|U\circ \exp_{x_0}\left(\frac\cdot\La\right)-h\right\|_{C^l(B_R)}<\ep\,,
\]
as claimed.

\bibliographystyle{amsplain}

\begin{thebibliography}{99}\frenchspacing


\bibitem{BR}
J. Bourgain, Z. Rudnick, On the nodal sets of toral eigenfunctions, Invent. Math. 185 (2011) 199--237.

\bibitem{CanHen}
Y. Canzani, B. Hanin, Scaling limit for the kernel of the spectral projector and remainder estimates in the pointwise Weyl law, Anal. PDE 8 (2015) 1707--1731.

\bibitem{CS19}
Y. Canzani, P. Sarnak, Topology and nesting of the zero set components of monochromatic random waves, Comm. Pure Appl. Math. 72 (2019) 343--374.

\bibitem{CC}
J. Cilleruelo, A. C\'ordoba, Lattice point on ellipses, Duke Math. J. 76 (1994) 741--750.

\bibitem{JEMS}
A. Enciso, D. Hartley, D. Peralta-Salas, A problem of Berry and knotted zeros zeros in the eigenfunctions of the harmonic oscillator, J. Eur. Math. Soc. 20 (2018) 301--314.

\bibitem{EPS13}
A.~Enciso, D.~Peralta-Salas, Submanifolds that are level sets of solutions to a second-order elliptic PDE, Adv. Math. 249 (2013) 204--249.


\bibitem{Acta}
A. Enciso, D. Peralta-Salas, Existence of knotted vortex tubes in
steady Euler flows, Acta Math. 214 (2015) 61--134.

\bibitem{APDE}
A. Enciso, D. Peralta-Salas, Bounded solutions to the Allen--Cahn equation with level sets of any compact topology, Anal. PDE 9 (2016) 1433--1446.

\bibitem{ENS}
A. Enciso, D. Peralta-Salas, F. Torres de Lizaur, Knotted structures
in high-energy Beltrami fields on the torus and the sphere,
Ann. Sci. \'Ec. Norm. Sup. 50 (2017) 995--1016.

\bibitem{proceedings}
A. Enciso, D. Peralta-Salas, F. Torres de Lizaur, High-energy eigenfunctions of the Laplacian on the torus and the sphere with nodal sets of complicated topology, Springer Proc. Math. \& Stat. 346 (2021) 245--261.

\bibitem{Alvaro}
A. Enciso, D. Peralta-Salas, \'A. Romaniega, Asymptotics for the nodal
components of non-identically distributed random monochromatic waves,
Int. Math. Res. Not. 2020, rnaa178, https://doi.org/10.1093/imrn/rnaa178.


\bibitem{Dima}
  D. Jakobson, Quantum limits on flat tori, Ann. of Math.  145 (1997)
  235--266.

\bibitem{Riviere}
H.  Hezari, G. Rivi\`ere, Equidistribution of toral eigenfunctions along hypersurfaces, Rev. Mat. Iberoam. 36 (2020) 435--454.

\bibitem{Hormander}
L. H\"ormander, The spectral function of an elliptic operator, Acta Math. 121 (1968) 193--218.

\bibitem{Hor15}
L.~H{\"o}rmander.
\newblock {\em The analysis of linear partial differential operators I},
\newblock  Springer, New York, 2015.

\bibitem{Iwaniec}
H. Iwaniec, {\em Topics in classical automorphic forms}, AMS,
Providence, 1997.

\bibitem{Lagace}
J.  Lagac\'e,
Eigenvalue optimisation on flat tori and lattice points in anisotropically expanding domains,
Canad. J. Math. 72 (2020) 967--987.

\bibitem{SW}
P. Sarnak, I. Wigman, Topologies of nodal sets of random-band limited functions, Comm. Pure Appl. Math. 72 (2019) 275--342.

\bibitem{Sevryuk}
M.B. Sevryuk, KAM-stable Hamiltonians, J. Dyn. Control Sys. 1 (1995) 351--366.

\bibitem{Uhlen}
K. Uhlenbeck, Generic properties of eigenfunctions, Amer. J. Math. 98 (1976) 1059--1078.


\end{thebibliography}

\end{document}